\documentclass[a4paper,12pt]{article}

\usepackage{graphicx}
\usepackage{amsmath,amsthm}
\usepackage{amsfonts,amssymb}
\usepackage{amscd}
\usepackage{color}
\usepackage[cp1251]{inputenc}
\usepackage[matrix,arrow,curve]{xy}
\usepackage{tikz}
\usetikzlibrary{arrows}

\newtheorem{lemma}{Lemma}[section]
\newtheorem{prop}[lemma]{Proposition}
\newtheorem{theorem}[lemma]{Theorem}
\newtheorem{coro}[lemma]{Corollary}

\newtheorem{Def}[lemma]{Definition}

\newtheorem{ex}{Example}

\newcommand{\HT}{{\rm HTrip}}
\newcommand{\QR}{{\rm RQuiv}}
\newcommand{\coQR}{{\rm coRQuiv}}

\newcommand{\HA}{{\rm HAlg}}

\newcommand{\Hom}{{\rm Hom}}

\newcommand{\kk}{{\bf k}}

\DeclareMathOperator{\Ker}{{Ker}}
\DeclareMathOperator{\Coker}{{Coker}}

\newcommand{\Id}{{\rm Id}}
\newcommand{\Ext}{{\rm Ext}}
\newcommand{\Tor}{{\rm Tor}}

\newcommand{\Ann}{{\rm Ann}}

\renewcommand{\Im}{{\rm Im\,}}

\newcommand{\HH}{{\mathcal H}}
\newcommand{\QQ}{{\mathcal Q}}
\newcommand{\OO}{{\mathcal O}}

\newcommand{\FF}{{\mathcal F}}
\newcommand{\GG}{{\mathcal G}}
\newcommand{\TT}{{\mathcal T}}

\newcommand{\bS}{{\bf S}}
\newcommand{\bL}{{\bf \Lambda}}

\newcommand{\vp}{\varphi}

\newcommand{\la}{\langle}
\newcommand{\ra}{\rangle}

\newcommand{\ot}{\otimes}

\renewcommand{\le}{\leqslant}
\renewcommand{\ge}{\geqslant}
\newcommand{\upper}[1]{\left\lceil#1\right\rceil}

\newenvironment{Ex}{\begin{ex}
\rm }{\flushright
$\Box$\end{ex}}

\newenvironment{Proof}[1][Proof.]{\begin{trivlist}
\item[\hskip \labelsep {\bfseries #1}]}{\flushright
$\Box$\end{trivlist}}

\numberwithin{equation}{section}
\begin{document}
\title{$s$-homogeneous algebras via $s$-homogeneous triples}
\author{E. Marcos, Y. Volkov}
\date{}
\maketitle
\begin{abstract}
To study $s$-homogeneous algebras, we introduce the category of quivers with $s$-homogeneous corelations and the category of $s$-homogeneous triples.
We show that both of these categories are equivalent to the category of $s$-homogeneous algebras. We prove some properties of the elements of $s$-homogeneous triples and give some consequences for $s$-Koszul algebras.
Then we discuss the relations between the $s$-Koszulity and the Hilbert series of $s$-homogeneous triples.
We give some application of the obtained results to $s$-homogeneous algebras with simple zero component. We describe all $s$-Koszul algebras with one relation recovering the result of Berger and all $s$-Koszul algebras with one dimensional $s$-th component.
We show that if the  $s$-th Veronese ring of an $s$-homogeneous algebra has two generators, then it has at least two relations.
Finally, we classify all $s$-homogeneous algebras with $s$-th Veronese rings $\kk\langle x,y\rangle/(xy,yx)$ and $\kk\langle x,y\rangle/(x^2,y^2)$. In particular, we show that all of these algebras are not $s$-Koszul while their $s$-homogeneous duals are $s$-Koszul.
\end{abstract}

\section{Introduction}

All the algebras under consideration are graded algebras of the form $\Lambda = \kk Q/I$, where the grading is induced by the path length and, in particular, $I$ is a homogeneous ideal.
The $\Ext$-algebra of the algebra $\Lambda$ is the graded algebra $\bigoplus_{i\geq 0} \Ext^i_{\Lambda}(\Lambda_0, \Lambda_0)$, where $\Lambda_0$, as usually, denotes the right $\Lambda$-module $\Lambda/\Lambda_{>0}$.

The notion of a Koszul algebra was introduced by S.~Priddy  in \cite{P}. All Koszul algebras are quadratic and they appear in 
pairs. The $\Ext$-algebra of a Koszul algebra $\Lambda$, is another time a Koszul algebra and its $\Ext$-algebra is the original
one. There is a duality between the category of Koszul modules over a Koszul algebra and the category of Koszul modules
over its $\Ext$-algebra. If we know that an algebra $\Gamma $ is the $\Ext$-algebra of a Koszul algebra $\Lambda$ then one may
recuperate $\Lambda$ from $\Gamma$ by taking its quadratic dual.
An interesting result states that an algebra is Koszul if and only if its $\Ext$-algebra is generated by its components of degrees 
zero and one. Note that even if a quadratic algebra is not Koszul, then it can be recovered from its quadratic dual algebra that in the non Koszul case is not anymore isomorphic to the $\Ext$-algebra.
Note also that if we have a quadratic algebra, then it is an $\Ext$-algebra of some quadratic algebra if and only if it is Koszul.

 Suppose now that the algebra $\Lambda = \kk Q/I$ is $s$-homogeneous, i.e. $I$ is generated by elements of degree $s$.
The notion of $s$-Koszul algebra was the first generalization of the notion of a Koszul algebra, and it was given for the first time in \cite{B} for a quiver with one vertex.
Later the definition was rewritten for the case of an arbitrary quiver in \cite{GMMZ}. It was shown in  the last mentioned work that if an algebra is $s$-homogeneous, then it is $s$-Koszul if and only if its $\Ext$-algebra is 
generated by its components of degree less than or equal to two.  The notion of $s$-Koszul algebra is important. For example, it was shown in \cite{BS} that an Artin-Shelter regular algebra of global
  dimension 3  is 3-Calabi-Yau if and only if it is $s$-Koszul.

The idea of the current paper appeared from the following question. Can we recover the algebra $\Lambda$ from its $\Ext$-algebra if $\Lambda$ is $s$-Koszul?
Note that, by the results of \cite{GMMZ}, the $\Ext$-algebra of $\Lambda$ is isomorphic in the $s$-Koszul case to the semidirect product of $(\Lambda^!)_0$ and $(\Lambda^!)_1$ after some regrading, where $\Lambda^!$ is the $s$-homogeneous dual of the algebra $\Lambda$.
In other words, $\Ext$-algebra of $\Lambda$ consists of the quadratic algebra $(\Lambda^!)_0$ and $(\Lambda^!)_0$-bimodule $(\Lambda^!)_1$. At this moment new interesting questions appear. Assume that we have a quadratic algebra $A$ and an $A$-bimodule $M$. Is  any $s$-homogeneous algebra $\Lambda$ such that the pair $(A,M)$ coincides with the pair $\big((\Lambda^!)_0,(\Lambda^!)_1\big)$? Can $\Lambda$ or $\Lambda^!$ be $s$-Koszul in this case? Should $\Lambda$ or $\Lambda^!$ be $s$-Koszul in this case?

Essentially, these questions appear, because the algebra structure of the $\Ext$ algebra does not contain all the $A_{\infty}$-structure in the $s$-Koszul case for $s>2$. Indeed, there is a missing map from $(\Lambda^!)_1^{\ot_{(\Lambda^!)_0}}$ to $(\Lambda^!)_0$ that is included in the algebra structure in the quadratic case. The notion of $s$-homogeneous triple comes from this argument and becomes a good alternative to the $s$-homogeneous dual in the $s$-homogeneous case. While an $s$-homogeneous algebra can be presented by a quiver with $s$-homogeneous relations, we will show that an $s$-homogeneous triple can be naturally presented by a quiver with $s$-homogeneous corelations. We will show that the pair $(A,M)$ has to satisfy some restrictive conditions to have a complement to an $s$-homogeneous triple. Moreover, we will show that in many cases the pair $\big((\Lambda^!)_0,(\Lambda^!)_1\big)$ determines the algebras $\Lambda$ and $\Lambda^!$.

The connections between the $s$-Koszulity and Hilbert series of an algebra was discussed in \cite{BS}. In particular, it was shown that Hilbert series of $s$-Koszul algebra and its Koszul dual satisfy some condition.
In this paper we give a further discussion of this. We rewrite the condition on Hilbert series of algebras in terms of Hilbert series of $s$-homogeneous triples and show that some part of this condition is equivalent to the extra condition introduced in \cite{B}.

In the last part of our work we show some applications of our technique. In our examples we consider only quivers with one vertex, i.e. algebras of the form $\Lambda=\kk\langle x_1,\dots,x_n\rangle/I$, where $I$ is an ideal generated by elements of degree $s$.
Firstly, we consider the case where $\dim_{\kk}I_s=1$. We discuss the $s$-Koszulity of $\Lambda$ in this case recovering and clarifying the results of \cite{B2} and show that $\Lambda^!$ is $s$-Koszul only in the case $n=1$ (i.e. in the case $\Lambda=\kk[x]/(x^s)$).
Next we consider the case $\dim_{\kk}I_s=2$. We show that in this case the $s$-th Veronese ring of $\Lambda^!$ is a quadratic algebra with two generators and at least two relations and classify all the algebras $\Lambda$, for which this ring is isomorphic to $\kk\langle x,y\rangle/(xy,yx)$ or $\kk\langle x,y\rangle/(x^2,y^2)$. It occurs that all such $\Lambda$ are $s$-Koszul, while $\Lambda^!$ is not $s$-Koszul in all cases.

\section{$s$-homogeneous relations and corelations}

We fix some notation during the paper. First of all, we fix some ground field $\kk$. Secondly, we fix some semisimple $\kk$-algebra $\OO$. Moreover, we assume that $\OO$ is isomorphic to a direct sum of copies of $\kk$ as an algebra and fix some basis $e_1,\dots,e_D$ of $\OO$ such that $e_ie_j=\delta_{i,j}e_j$. Everything in this paper is over $\OO$. So we write simply $\otimes$ instead of $\otimes_\OO$. All modules in this paper are right modules if the opposite is not stated. If $U$ is an $\OO$-bimodule, then $T(U)=\oplus_{i\ge 0}U^{\otimes i}$ denotes the tensor algebra of $U$ over $\OO$. For convenience we set everywhere $U^{\ot i}=0$ for $i<0$. If $A$ is an algebra, $M$ is a right $A$-module, and $X$ is a subset of $M$, then $\langle X\rangle_A$ denotes the right $A$-submodule of $M$ generated by the set $X$.

\begin{Def} {\rm An {\it $\OO$-quiver with relations} is a triple $(U,V,\iota)$, where $U$ and $V$ are $\OO$-bimodules and $\iota:V\hookrightarrow \oplus_{i\ge 2}U^{\ot i}$ is an $A$-bimodule monomorphism. In this situation $U$ is called an {\it $\OO$-quiver} and $\Im\iota$ is called a {\it set of relations}. The $\OO$-quiver $U$ is called {\it finite} if $U$ is a finitely generated $\OO$-bimodule. The set of relations $\Im\iota$ is called {\it $s$-homogeneous} if $\Im\iota\subset U^{\ot s}$.  If $(U',V',\iota')$ is another $\OO$-quiver with relations, then a {\it morphism} from $(U,V,\iota)$ to $(U',V',\iota')$ is a pair $(f,g)$, where $f:U\rightarrow U'$ and $g:V\rightarrow V'$ are such $\OO$-bimodule homomorphisms that $\iota'g=(\oplus_{i\ge 2}f^{\ot i})\iota$. We denote by $\QR(\OO,s)$ the category of finite $\OO$-quivers with {\it $s$-homogeneous} sets of relations.
}
\end{Def}

Since $\OO$ is fixed, we will write simply quiver instead of $\OO$-quiver. Now we are going to define a quiver with corelations. Since the present paper is devoted to $s$-homogeneous case, we define only $s$-homogeneous corelations.

\begin{Def} {\rm A {\it quiver with $s$-homogeneous corelations} is a triple $(U,W,\pi)$, where $U$ and $W$ are $\OO$-bimodules and $\pi:U^{\ot s}\twoheadrightarrow W$ is an $A$-bimodule epimorphism. In this situation $U$ is called a {\it quiver} as before and $W$ is called a {\it set of corelations}.  If $(U',W',\pi')$ is another quiver with $s$-homogeneous corelations, then a {\it morphism} from $(U,W,\pi)$ to $(U',W',\pi')$ is a pair $(f,g)$, where $f:U\rightarrow U'$ and $g:W\rightarrow W'$ are such $\OO$-bimodule homomorphisms that 
$g\pi=\pi' f^{\ot s}$. We denote by $\coQR(\OO,s)$ the category of finite quivers with {\it $s$-homogeneous} sets of corelations.
}
\end{Def}

Note that $\mathcal{D}=\Hom_{\kk}(-,\kk)$ is a contravariant endofunctor of the category of finitely generated $\OO$-bimodules. Moreover, $\mathcal{D}^2$ is isomorphic to the identity functor, i.e. $\mathcal{D}$ is a duality. We will write $(-)^*$ instead of $\mathcal{D}(-)$. Note that, for finitely generated $\OO$-bimodules $M$ and $L$, $(M\ot L)^*$ can be identified with $M^*\ot L^*$ via the canonical embedding $M^*\ot L^*\hookrightarrow (M\ot L)^*$. Thus, we obtain also the contravariant functor $\mathfrak{D}_R:\QR(\OO,s)\rightarrow \coQR(\OO,s)$ defined by the equalities $\mathfrak{D}_R(U,V,\iota)=(U^*,V^*,\iota^*)$ and $\mathfrak{D}_R(f,g)=(f^*,g^*)$ for an object $(U,V,\iota)$ and a morphism $(f,g)$ of the category $\QR(\OO,s)$. Analogously one can construct $\mathfrak{D}_{coR}:\coQR(\OO,s)\rightarrow \QR(\OO,s)$. We will write $(-)^*$ instead of $\mathfrak{D}_R(-)$ and $\mathfrak{D}_{coR}(-)$.

Now we define also the functors $\mathfrak{Coker}:\QR(\OO,s)\rightarrow \coQR(\OO,s)$ and $\mathfrak{Ker}:\coQR(\OO,s)\rightarrow \QR(\OO,s)$ in the following way. We simply define $\mathfrak{Coker}(U,V,\iota)=(U,\Coker\iota,\pi)$, where $\pi:U^{\ot s}\twoheadrightarrow\Coker\iota$ is the canonical projection and $\mathfrak{Ker}(U,W,\pi)=(U,\Ker\pi,\iota)$, where $\iota:\Ker\pi\hookrightarrow U^{\ot s}$ is the canonical inclusion. The definition on morphisms is the natural one, we simply define $\mathfrak{Coker}(f,g)=\big(f,\Coker(g,f^{\ot s})\big)$ and $\mathfrak{Ker}(f,g)=\big(f,\Ker(f^{\ot s},g)\big)$ for a morphism $(f,g)$ in the corresponding category. We have the following lemma.

\begin{lemma}\label{DandK} $(\mathfrak{D}_R, \mathfrak{D}_{coR})$ is a pair of quasi inverse dualities and $(\mathfrak{Coker},\mathfrak{Ker})$ is a pair of quasi inverse equivalences. Moreover, there is an isomorphism of contravariant functors $\mathfrak{Ker}\circ\mathfrak{D}_R\cong \mathfrak{D}_{coR}\circ\mathfrak{Coker}$.
\end{lemma}

We will write $(-)^!$ instead of $\mathfrak{Ker}\circ\mathfrak{D}_R(-)$ and $\mathfrak{Coker}\circ\mathfrak{D}_{coR}(-)$.

\section{$s$-homogeneous algebras}\label{sAlg}

In this section we introduce the notion of $s$-homogeneous algebra. Usually, it is defined as an algebra corresponding to an element of $\QR(\OO,s)$ as will be explained, but we prefer to give an intrinsic definition that does not use any presentation of an algebra.

\begin{Def}{\rm An {\it $\OO$-algebra} is a graded algebra $\Lambda=\oplus_{i\ge 0}\Lambda_i$ such that $\Lambda_0=\OO$.
Given an $\OO$-algebra $\Lambda$, we denote by $\Lambda^{(r)}$ the {\it $r$-Veronese ring} of $\Lambda$, i.e. a graded algebra $\Lambda^{(r)}=\oplus_{i\ge 0}\Lambda^{(r)}_i$, where $\Lambda^{(r)}_i=\Lambda_{ri}$ and the multiplication of $\Lambda^{(r)}$ is induced by the multiplication of $\Lambda$.
Also we define the {\it $(r,t)$-Veronese bimodule} $\Lambda^{(r,t)}$ as a graded $\Lambda^{(r)}$-bimodule $\Lambda^{(r,t)}=\oplus_{i\ge 0}\Lambda^{(r,t)}_i$, where $\Lambda^{(r,t)}_i=\Lambda_{ri+t}$ and the $\Lambda^{(r)}$-bimodule structure on $\Lambda^{(r,t)}$ is induced by the multiplication of $\Lambda$.
}
\end{Def}

Note that the multiplication of $\Lambda$ induces a $\Lambda^{(r)}$-bimodule homomorphism from $\left(\Lambda^{(r,1)}\right)^{\ot^{r}_{\Lambda^{(r)}}}$ to $\Lambda^{(r)}$. We denote this homomorphism by $\phi_{\Lambda}^{(r)}$. Also we denote by $\phi_{\Lambda}^{r}$ the zero component of $\phi_{\Lambda}^{(r)}$, that is a map from $\Lambda_1^{\ot r}$ to $\Lambda_r$. A {\it morphism} of $\OO$-algebras is by definition a graded homomorphism of algebras identical on the zero component $\OO$. Any such a morphism induces homomorphisms between all Veronese rings and all Veronese bimodules of the algebras under consideration. Moreover, the induced morphisms are compatible with the maps $\phi_{\Lambda}^{(r)}$.

\begin{Def}{\rm Given an $\OO$-algebra $\Lambda$, we call it {\it $s$-homogeneous} if $\phi_{\Lambda}^r$ is surjective and
$$\Ker\phi_{\Lambda}^r=\sum\limits_{i=0}^{r-s}\Lambda_1^{\ot i}\ot\Ker\phi_{\Lambda}^s\ot\Lambda_1^{\ot (r-s-i)}$$
for any integer $r\ge 1$. In particular, if $\Lambda$ is $s$-homogeneous, then $\phi_{\Lambda}^r$ is bijective for $r<s$. We denote by $\HA(\OO,s)$ the category of $s$-homogeneous $\OO$-algebras.
We will write sometimes quadratic algebras instead of $2$-homogeneous algebras.
}
\end{Def}

There is another definition of an $s$-homogeneous algebra in terms of its Veronese modules. We give it in our next statement.

\begin{prop}\label{shom} Let $\Lambda$ be an $\OO$-algebra, and, for $n,r,t\ge 1$,
$$\phi_n^{r,t}:\Lambda^{(n,r)}\ot_{\Lambda^{(n)}}\Lambda^{(n,t)}\rightarrow \Lambda^{(n,r+t)}$$ be the $\Lambda^{(n)}$-bimodule  homomorphism induced by the multiplication of $\Lambda$. Then the following conditions are equivalent:
\begin{enumerate}
\item $\Lambda\in\HA(\OO,s)$.
\item $\phi_n^{r,t}$ is surjective for all $n,r,t\ge 1$; $\phi_n^{r,t}$ is bijective if $n\ge s-1$ and $r+t<s$; $\phi_n^{r,t}$ splits and $\Ker\phi_n^{r,t}$ is concentrated in degree $0$ if $n\ge s-1$ and $r+t\ge s$.
\item $(\phi_s^{r,1})_0$ is surjective for any $r\ge 1$ and is bijective for $r<s-1$; $(\phi_n^{1,1})_1$ is bijective for $n\ge s-1$.
\end{enumerate}
\end{prop}
\begin{Proof} For any integer $i,j,k\ge 0$ such that $j+k=i$ and $n,r,t\ge 1$, we have a map $$\mu^{n,j,k}_{r,t}:{\Lambda}_{nj+r}\ot {\Lambda}_{nk+t}\rightarrow (\Lambda^{(n,r)}\ot_{\Lambda^{(n)}}\Lambda^{(n,t)})_i.$$
Then we have the commutative diagram
\begin{center}
\begin{tikzpicture}
 \node(L1) {${\Lambda}_1^{\ot (ni+r+t)}$};
\node(L11) [below  of=L1]{};
 \node(jxk) [below  of=L11]{${\Lambda}_{ni+r}\ot {\Lambda}_{t}$};
\node(jxk1) [right of=jxk]{};
\node(jxk2) [right of=jxk1]{};
\node(jxk3) [right of=jxk2]{};
\node(jxk4) [right of=jxk3]{};
\node(jxk5) [right of=jxk4]{};
 \node(rxt) [right of=jxk5]{$(\Lambda^{(n,r)}\ot_{\Lambda^{(n)}}\Lambda^{(n,t)})_i$};
\node(L111) [right of=L1]{};
\node(L112) [right of=L111]{};
\node(L113) [right of=L112]{};
\node(L114) [right of=L113]{};
\node(L115) [right of=L114]{};
\node(j+k) [right of=L115]{${\Lambda}_{ni+r+t}$};
\node(j+k1) [right of=j+k]{};
\node(j+k2) [right of=j+k1]{};
\node(j+k3) [right of=j+k2]{};
\node(j+k4) [right of=j+k3]{};
\node(r+t) [below of=j+k4]{$\Lambda^{(n,r+t)}_i$};

\draw [double equal sign distance] (j+k)  to (r+t);

\draw [->,>=stealth'] (L1)   to node[left, fill=white]{\tiny$\phi_{\Lambda}^{ni+r}\ot \phi_{\Lambda}^{t}$ } (jxk) ;equ

\draw [->,>=stealth'] (L1)   to node[above, fill=white]{\tiny$\phi_{\Lambda}^{ni+r+t}$ } (j+k) ;

\draw [->,>=stealth'] (rxt)   to node[above=7, left=3, fill=white]{\tiny$(\phi_n^{r,t})_i$ } (r+t) ;

\draw [->,>=stealth'] (jxk)  to node[above, fill=white]{\tiny$\mu^{n,i,0}_{r,t}$ } (rxt) ; 
\end{tikzpicture}
\end{center}

``$1\implies 2$" The surjectivity of $(\phi_n^{r,t})_i$ follows from the surjectivity of $\phi_{\Lambda}^{ni+r+t}$ and the just mentioned commutative diagram.

Suppose that $j>0$, $k\ge 0$. It is clear that $\Lambda_1^{\ot(ni+r+t)}$ is generated as a $\kk$-linear space by the elements of the form $a\ot b \ot c$, where $a\in \Lambda_1^{\ot (n(j-1)+r)}$, $b\in\Lambda_1^{\ot n}$ and $c\in\Lambda_1^{\ot (nk+t)}$. But for such an element we have
\begin{multline*}
\mu^{n,j,k}_{r,t}(\phi_{\Lambda}^{nj+r}(a\ot b)\ot \phi_{\Lambda}^{nk+t}(c))=\mu^{n,j,k}_{r,t}(\phi_{\Lambda}^{n(j-1)+r}(a)\phi_{\Lambda}^{n}(b)\ot \phi_{\Lambda}^{nk+t}(c))\\
=\mu^{n,j-1,k+1}_{r,t}(\phi_{\Lambda}^{n(j-1)+r}(a)\ot \phi_{\Lambda}^{n(k+1)+t}(b\otimes c)).
\end{multline*}
Thus, we have $\mu^{n,j,k}_{r,t}(\phi_{\Lambda}^{nj+r}\ot \phi_{\Lambda}^{nk+t})=\mu^{n,i,0}_{r,t}(\phi_{\Lambda}^{ni+r}\ot \phi_{\Lambda}^{t})$.
Since $\phi_{\Lambda}^{nj+r}\ot \phi_{\Lambda}^{nk+t}$ is surjective, this means, in particular, that $\Im\mu^{n,j,k}_{r,t}$ depends only on the sum of $j$ and $k$, but does not depend on $j$ and $k$.
Thus, we have $(\Lambda^{(n,r)}\ot_{\Lambda^{(n)}}\Lambda^{(n,t)})_i=\sum_{j'+k'=i}\Im\mu^{n,j',k'}_{r,t}=\Im\mu^{n,j,k}_{r,t}$. In particular, $\mu^{n,j,k}_{r,t}$ is surjective for any values of indices.

From now we assume that $n\ge s-1$. Suppose that $(\phi_n^{r,t})_i(u)=0$ for some $u\in(\Lambda^{(n,r)}\ot_{\Lambda^{(n)}}\Lambda^{(n,t)})_i$. Since $\phi_{\Lambda}^{ni+r}\ot \phi_{\Lambda}^{t}$ and $\mu^{n,i,0}_{r,t}$ are surjective, we have $u=\mu^{n,i,0}_{r,t}(\phi_{\Lambda}^{ni+r}\ot \phi_{\Lambda}^{t})(v)$ for some $v\in{\Lambda}_1^{\ot (ni+r+t)}$. It follows from the commutative diagram above that $\phi_{\Lambda}^{ni+r+t}(v)=0$. Thus,
$$
v\in\Ker\phi_{\Lambda}^{ni+r+t}=\sum\limits_{l=0}^{ni-s+r+t}\Lambda_1^{\ot l}\ot\Ker\phi_{\Lambda}^s\ot\Lambda_1^{\ot (ni-s+r+t-l)}.
$$
Let us prove that $\mu^{n,i,0}_{r,t}(\phi_{\Lambda}^{ni+r}\ot \phi_{\Lambda}^{t})(\Lambda_1^{\ot l}\ot\Ker\phi_{\Lambda}^s\ot\Lambda_1^{\ot (ni-s+r+t-l)})=0$ if $i>0$. Let $m$ denote $ni-s+r$. If $l\le m$, then we have
$$
(\phi_{\Lambda}^{m+s}\ot \phi_{\Lambda}^{t})(\Lambda_1^{\ot l}\ot\Ker\phi_{\Lambda}^s\ot\Lambda_1^{\ot (m+t-l)})=\phi_{\Lambda}^l(\Lambda_1^{\ot l})\phi_{\Lambda}^s(\Ker\phi_{\Lambda}^s)\phi_{\Lambda}^{m-l}(\Lambda_1^{\ot (m-l)})\ot \phi_{\Lambda}^{t}(\Lambda_1^{\ot t})=0.
$$
If $l\ge m+1$, then we have $\mu^{n,i,0}_{r,t}(\phi_{\Lambda}^{m+s}\ot \phi_{\Lambda}^{t})=\mu^{n,i-1,1}_{r,t}(\phi_{\Lambda}^{m+s-n}\ot \phi_{\Lambda}^{n+t})$ and
\begin{multline*}
(\phi_{\Lambda}^{m+s-n}\ot \phi_{\Lambda}^{n+t})(\Lambda_1^{\ot l}\ot\Ker\phi_{\Lambda}^s\ot\Lambda_1^{\ot (m+t-l)})\\=\phi_{\Lambda}^{m+s-n}(\Lambda_1^{\ot (m+s-n)})\ot\phi_{\Lambda}^{l-m-s+n}(\Lambda_1^{\ot (l-m-s+n)}) \phi_{\Lambda}^s(\Ker\phi_{\Lambda}^s)
\phi_{\Lambda}^{m+t-l}(\Lambda_1^{\ot (m+t-l)})=0.
\end{multline*}
Thus, $u=\mu^{n,i,0}_{r,t}(\phi_{\Lambda}^{ni+r}\ot \phi_{\Lambda}^{t})(v)=0$ and $(\phi_n^{r,t})_i$ is bijective if $i>0$.

If $r+t<s$, then $\Ker\phi_{\Lambda}^{r+t}=0$, i.e. $(\phi_n^{r,t})_0$ is also bijective.

It remains to prove that $\phi_n^{r,t}$ splits if $r+t\ge s$. Let us denote $\Lambda^{(n,r)}\ot_{\Lambda^{(n)}}\Lambda^{(n,t)}$ by $M$. It is clear that $M$ a graded $\Lambda^{(n)}$-bimodule concentrated in nonnegative degrees. Let $\pi:M\twoheadrightarrow M/M_{>0}$ be the canonical projection and $\iota:\Ker\phi_n^{r,t}\hookrightarrow M$ be the canonical inclusion. It is enough to show that $\iota$ splits. Since $\Ker\phi_n^{r,t}$ is concentrated in degree $0$, the map $\pi\iota$ is a monomorphism. Since $M/M_{>0}$ is semisimple, $\pi\iota$ splits. Consequently, $\iota$ splits too.

``$2\implies 3$" is clear.

``$3\implies 1$" Note that $\mu^{n,0,0}_{r,t}:{\Lambda}_{r}\ot {\Lambda}_{t}\rightarrow (\Lambda^{(n,r)}\ot_{\Lambda^{(n)}}\Lambda^{(n,t)})_0$ is an isomorphism for any $n,r,t\ge 1$.

Using the commutative diagram above and the fact that $(\phi_s^{r,1})_0$ is surjective for any $r\ge 1$, we get surjectivity of $\phi_{\Lambda}^n$ for any $n\ge 1$ by induction on $n$.
Using the commutative diagram above and the fact that $(\phi_s^{r,1})_0$ is bijective for any $r<s-1$, we get bijectivity of $\phi_{\Lambda}^n$ for any $n<s$ by induction on $n$.

Now it is enough to show that $\Ker\phi_{\Lambda}^n=\Ker\phi_{\Lambda}^{n-1}\otimes \Lambda_1+\Lambda_1\otimes \Ker\phi_{\Lambda}^{n-1}$ for $n>s$. Let us take some $n>s$ and $u\in\Ker\phi_{\Lambda}^n$.
We have $(\phi_{n-2}^{1,1})_1\mu_{1,1}^{n-2,1,0}(\phi_{\Lambda}^{n-1}\otimes \phi_{\Lambda}^{1})(u)=0$. Since $(\phi_{n-2}^{1,1})_1$ is bijective, we have $\mu_{1,1}^{n-2,1,0}(\phi_{\Lambda}^{n-1}\otimes \phi_{\Lambda}^{1})(u)=0$.
It is easy to see that the kernel of the canonical map from $(\Lambda_{n-1}\ot \Lambda_1)\oplus(\Lambda_{1}\ot \Lambda_{n-1})$ to $(\Lambda^{(n-2,1)}\ot_{\Lambda^{(n-2)}}\Lambda^{(n-2,1)})_1$ is
$(\phi_{\Lambda}^{n-1}\otimes \phi_{\Lambda}^{1}-\phi_{\Lambda}^{1}\otimes \phi_{\Lambda}^{n-1})(\Lambda_1^{\ot n})$. Thus, $(\phi_{\Lambda}^{n-1}\otimes \phi_{\Lambda}^{1})(u)=(\phi_{\Lambda}^{n-1}\otimes \phi_{\Lambda}^{1})(u')$ for some $u'\in\Lambda_1\otimes \Ker\phi_{\Lambda}^{n-1}$ and we have $u=(u-u')+u'\in\Ker\phi_{\Lambda}^{n-1}\otimes \Lambda_1+\Lambda_1\otimes \Ker\phi_{\Lambda}^{n-1}$.
\end{Proof}

Let us define the {\it representing functor} $\mathfrak{ARep}:\HA(\OO,s)\rightarrow\QR(\OO,s)$ and the {\it algebralizing functor} $\mathfrak{Alg}:\QR(\OO,s)\rightarrow\HA(\OO,s)$ in the following way. Given $\Lambda\in\HA(\OO,s)$, we define $\mathfrak{ARep}(\Lambda)=(\Lambda_1,\Ker\phi_{\Lambda}^s,\iota)$, where $\iota:\Ker\phi_{\Lambda}^s\hookrightarrow \Lambda_1^{\ot s}$ is the canonical inclusion. If $\Lambda'$ is another $s$-homogeneous $\OO$-algebra and $f:\Lambda\rightarrow\Lambda'$ is a morphism of $\OO$-algebras, then we define $\mathfrak{ARep}(f)=(f|_{\Lambda_1},(f|_{\Lambda_1})^{\ot s}|_{\Ker\phi_{\Lambda}^s})$.

Let us now consider $(U,V,\iota)\in\QR(\OO,s)$. We define $\mathfrak{Alg}(U,V,\iota)=T(U)/I_{\iota}$, where $I_{\iota}$ is the ideal of $T(U)$ generated by $\Im\iota\subset U^{\ot s}$.
It is clear that $T(U)/I_{\iota}\in\HA(\OO,s)$. If $(f,g):(U,V,\iota)\rightarrow (U',V',\iota')$ is a morphism of quivers with relations, then it induces the morphism $T(f):T(U)\rightarrow T(U')$ such that $T(f)|_{U^{\ot i}}=f^{\ot i}$ for all $i\ge 0$. It is easy to see that $T(f)(I_{\iota})\subset I_{\iota'}$, i.e. $T(f)$ induces a well defined morphism of $\OO$-algebras $$\mathfrak{Alg}(f,g):\mathfrak{Alg}(U,V,\iota)=T(U)/I_{\iota}\rightarrow T(U')/I_{\iota'}=\mathfrak{Alg}(U',V',\iota').$$
The proof of the next proposition is standard.

\begin{prop}\label{RA} $(\mathfrak{ARep},\mathfrak{Alg})$ is a pair of quasi inverse equivalences.
\end{prop}

Thus, we have a duality $(-)^!=\mathfrak{Alg}\circ (-)^!\circ\mathfrak{ARep}:\HA(\OO, s)\rightarrow \HA(\OO, s)$. Given $\Lambda\in\HA(\OO,s)$, we will call the algebra $\Lambda^!$ the {\it $s$-homogeneous dual algebra} for $\Lambda$ or {\it $s$-dual algebra} for $\Lambda$ for short. Note that the definition depends on $s$, i.e. if $\Lambda\in\HA(\OO,s)$ and $\Lambda\in\HA(\OO,s')$, then the $s$-dual and $s'$-dual algebras of $\Lambda$ are not isomorphic. This will not cause any confusion, because we always fix $s$.

\section{$s$-homogeneous triples}

Another object that we are going to study in this paper is the $s$-homogeneous triples. In this section we introduce their definition and show that they can be naturally represented by quivers with corelations.

Given an $\OO$-algebra $A$, the category of graded $A$-modules is the category whose objects are graded $A$-modules and whose morphisms are degree preserving homomorphisms of $A$-modules.
A graded $A$-module $M$ is called {\it linear until the $n$-th degree} if there exists a projective resolution of $M$ in the category of graded $A$-modules
$$
M\leftarrow P_0\leftarrow P_1\leftarrow\cdots\leftarrow P_n\leftarrow\cdots
$$
such that $P_i$ is generated in degree $i$, i.e. $P_i=(P_i)_iA$, for $0\le i\le n$.

\begin{Def}\label{str}
{\rm
An {\it $s$-homogeneous  triple} is a triple $(A,M,\vp)$, where $A$ is a quadratic $\OO$-algebra, $M$ is a graded $A$-bimodule which is linear until the first degree as left and as right $A$-module, and $\vp:M^{\ot^{s}_A }\rightarrow A(1)$ is a homomorphism of graded $A$-bimodules such that 
\begin{enumerate}
\item $\Im\vp=A_{>0}(1)$;
\item $1_M\ot_A\vp=\vp\ot_A 1_M:M^{\ot_A^{s+1}}\rightarrow M(1)$;
\item $\Ker(1_M\ot_A\vp)=\Ker\vp\ot_A M+M\ot_A\Ker\vp$;
\item $\Ker(\vp\ot_A\vp)=\sum\limits_{i=0}^sM^{\ot_A^i}\ot_A\Ker\vp\ot_A M^{\ot_A^{s-i}}$.
\end{enumerate}
}
\end{Def}

Let $(A,M,\vp)$ and $(B,L,\psi)$ be $s$-homogeneous  triples. A {\it morphism} from $(A,M,\vp)$ to $(B,L,\psi)$ is a pair $(f,g)$, where $f:A\rightarrow B$ is a morphism of graded $\OO$-algebras and $g:M\rightarrow L$ is a morphism of graded $A$-bimodules such that $f\vp=\psi g^{\ot s}$. Here the $A$-bimodule structure on $L$ is induced by the map $f$. Let $\HT(s,\OO)$ denote the category of $s$-homogeneous triples.

Note that if $(A,M,\vp)$ is an $s$-homogeneous triple and $n\ge s$ is an integer, then the map $1_{M^{\ot_A^i}}\ot_A\vp\ot_A1_{M^{\ot_A^{n-s-i}}}:M^{\ot_A^n}\rightarrow M^{\ot_A^{n-s}}(1)$ does not depend on $0\le i\le n-s$ due to the second point of Definition \ref{str}. For simplicity we denote this map by $\vp$ too. So, if $n\ge ks$, then we have a graded $A$-bimodule homomorphism $\vp^k:M^{\ot_A^n}\rightarrow M^{\ot_A^{n-ks}}(s)$. As usually, $M^{\ot_A^0}=A$ everywhere.

Let us define the {\it representing functor} $\mathfrak{TRep}:\HT(\OO,s)\rightarrow\coQR(\OO,s)$ and the {\it triplizing functor} $\mathfrak{Trip}:\coQR(\OO,s)\rightarrow\HT(\OO,s)$ in the following way. Given $(A,M,\vp)\in\HT(\OO,s)$, we define $\mathfrak{TRep}(A,M,\vp)=(M_0,A_1,\vp_0)$. If $(A',M',\vp')$ is another $s$-homogeneous triple and $(f,g):(A,M,\vp)\rightarrow(A',M',\vp')$ is a morphism in $\HT(\OO,s)$, then we define $\mathfrak{TRep}(f,g)=(g_0,f_1)$.


Let us now consider the quiver with corelations $\QQ=(U,W,\pi)$. Let us define $A_{\QQ}=T(W)/I_{\pi}$, where $I_{\pi}$ is the ideal of $T(W)$ generated by $$(\pi\ot \pi)\left(\sum\limits_{i=0}^sU^{\ot i}\ot\Ker\pi\ot U^{\ot(s-i)}\right)\subset W\ot W.$$
Let us equip all the elements of $U$ with degree $0$ grading and consider the graded projective $A_{\QQ}$-bimodule $P_{\QQ}=A_{\QQ}\ot U\ot A_{\QQ}$. Since $(A_{\QQ})_1=W$, the sets $W\ot U$ and $U\ot W$ can be naturally identified with subspaces $(A_{\QQ})_1\ot U\ot \OO$ and $\OO\ot U\ot (A_{\QQ})_1$ of $P_{\QQ}$. Let us define $M_{\QQ}=P_{\QQ}/N_{\QQ}$, where $N_{\QQ}=\left\langle (\pi\ot 1_U-1_U\ot\pi)(U^{\ot(s+1)})\right\rangle_{A^{\rm op}\ot A}\subset P_{\QQ}$. For any $m\ge 1$ let us consider the projective right $A_{\QQ}$-module $P_m^R=U^{\ot m}\ot A_{\QQ}$.
There is a right $A_{\QQ}$-module homomorphism $\beta_m:P_m^R\rightarrow M_{\QQ}^{\ot_{A_{\QQ}}^m}$ that sends $u_1\ot \cdots \ot u_m\ot a\in P_m^R$ to the class of the element
$$(1_{\OO}\ot u_1\ot 1_{\OO})\ot_{A_{\QQ}}\cdots\ot_{A_{\QQ}}(1_{\OO}\ot u_{m-1}\ot 1_{\OO})\ot_{A_{\QQ}}(1_{\OO}\ot u_{m}\ot a).$$
It easily follows from the definition of $N_{\QQ}$ that $\beta_m$ is surjective for any $m\ge 1$. Let us denote by $\rho_i:U^{\ot (is)}\rightarrow (A_{\QQ})_i$ the composition of $\pi^{\ot i}$ and the natural projection $W^{\ot i}\twoheadrightarrow (A_{\QQ})_i$.

\begin{lemma}\label{right}
The kernel of the map $\beta_1$ is $\left\langle (1_U\ot\pi)(\Ker\pi\ot U)\right\rangle_A$. In particular, $M_{\QQ}$ is linear until the first degree as a right $A_{\QQ}$-module.
\end{lemma}
\begin{Proof} Let $M$ denote the right $A_{\QQ}$-module $(U\ot A_{\QQ})/\left\langle (1_U\ot\pi)(\Ker\pi\ot U)\right\rangle_A$. Direct verification shows that $\beta_1$ induces an epimorphism $\beta:M\rightarrow M_{\QQ}$.
We define $\alpha:M_{\QQ}\rightarrow M$ in the following way. Let us consider the map from $P_{\QQ}$ to $M$ that sends the element of the form $a\ot u\ot b\in (A_{\QQ})_i\ot U\ot A_{\QQ}$, to the class of the element $(1_U\ot\rho_i)(v\ot u)a$, where $v\in U^{\ot (is)}$ is some element such that $\rho_i(v)=b$. Since two elements $v,v'\in U^{\ot (is)}$ such that $\rho_i(v)=\rho_i(v')=b$ differ by an element from
$\sum\limits_{j=0}^{is-s}U^{\ot i}\ot\Ker\pi\ot U^{\ot(is-s-j)}$, it is easy to see that $(1_U\ot\rho_i)(v\ot u-v'\ot u)$ belongs to $(1_U\ot \pi)(\Ker\pi\ot U)A_{\QQ}$. Thus, the map just defined  is a well defined homomorphism of right $A_{\QQ}$-modules. It is not difficult to see also that this map vanishes on $(\pi\ot 1_U-1_U\ot\pi)(U^{\ot(s+1)})$, and hence induces the required homomorphism $\alpha$. It is not difficult also  to see that $\alpha\beta=\Id_M$, and hence $\beta$ is an isomorphism.
\end{Proof}

The dual argument shows that $M_{\QQ}$ is linear until the first degree as a left $A_{\QQ}$-module.
Let us define the homomorphism of graded bimodules $\vp_{\QQ}:M_{\QQ}^{\ot^{s}_A}\rightarrow A(1)$ in the following way. For a homogeneous element
$$(a_1\ot u_1\ot b_1)\ot_A\cdots\ot_A(a_s\ot u_s\ot b_s)\in (A_{\QQ})_{i_1}\ot U\ot (A_{\QQ})_{j_1}\ot_A\cdots\ot_A(A_{\QQ})_{i_s}\ot U\ot (A_{\QQ})_{j_s} $$
 we choose elements $v_1,\dots,v_{s-1}$ in such a way that $\rho_{j_{k}+i_{k+1}}(v_k)=b_{k}a_{k+1}$ for $1\le k\le s-1$ and set
$$\tilde\vp_{\QQ}\big((a_1\ot u_1\ot b_1)\ot_A\cdots\ot_A(a_s\ot u_s\ot b_s)\big)=a_1\rho_{\sum\limits_{k=1}^{s-1}(j_{k}+i_{k+1})+1}(u_1\ot v_1\ot\cdots\ot u_{s-1}\ot v_{s-1}\ot u_s)b_s.$$
Analogously to the proof of Lemma \ref{right}, one can  show that $\tilde\vp_{\QQ}:P_{\QQ}^{\ot^{s}_A}\rightarrow A(1)$ is a well defined homomorphism of graded $A_{\QQ}$-bimodules. Moreover, it is clear that $\Im\tilde\vp_{\QQ}=A_{>0}(1)$.
Let $\iota_{\QQ}:N_{\QQ}\hookrightarrow P_{\QQ}$ be the canonical inclusion.
Since the kernel of the projection $P_{\QQ}^{\ot^{s}_A}\twoheadrightarrow M_{\QQ}^{\ot^{s}_A}$ belongs to the image of
$$\sum_{i=1}^s1_{P_{\QQ}}^{\ot^{i-1}_A}\ot_A\iota_{\QQ}\ot_A1_{P_{\QQ}}^{\ot^{s-i}_A}:\sum_{i=1}^sP_{\QQ}^{\ot^{i-1}_A}\ot_AN_{\QQ}\ot_AP_{\QQ}^{\ot^{s-i}_A}\rightarrow P_{\QQ}^{\ot^{s}_A}$$
    and it is easy to check that $\tilde\vp_{\QQ}\left(\sum_{i=1}^s1_{P_{\QQ}}^{\ot^{i-1}_A}\ot_A\iota_{\QQ}\ot_A1_{P_{\QQ}}^{\ot^{s-i}_A}\right)=0$, we get the homomorphism $\vp_{\QQ}$.

\begin{lemma} $(A_{\QQ},M_{\QQ},\vp_{\QQ})$ is an $s$-homogeneous triple.
\end{lemma}
\begin{Proof} It remains to prove Conditions 2-4 from the definition of an $s$-homogeneous triple. It is not difficult to verify Condition 2. In particular, we have
\begin{multline*}
\Ker\vp_{\QQ}\ot_{A_{\QQ}} M_{\QQ}+M_{\QQ}\ot_{A_{\QQ}}\Ker\vp_{\QQ}\subset\Ker(1_{M_{\QQ}}\ot_{A_{\QQ}}\vp_{\QQ})\mbox{ and}\\
\sum\limits_{i=0}^sM_{\QQ}^{\ot^{i}_{A_{\QQ}}}\ot_{A_{\QQ}}\Ker\vp_{\QQ}\ot_{A_{\QQ}} M_{\QQ}^{\ot^{s-i}_{A_{\QQ}}}\subset\Ker(\vp_{\QQ}\ot_{A_{\QQ}}\vp_{\QQ}).
\end{multline*}
Let us prove the reverse inclusions.
Suppose that $x\in M_{\QQ}^{\ot^{s+1}_{A_{\QQ}}}$ is such that $(1_{M_{\QQ}}\ot_{A_{\QQ}}\vp_{\QQ})(x)=0$. We may assume that $x$ is a homogeneous element of degree $k$.
Let us represent $x$ in the form $x=\beta_{s+1}(1_U^{\ot(s+1)}\ot \rho_k)(y)$ for $y\in U^{\ot (s+1+ks)}$.
We have
$$y\in\sum\limits_{i=0}^{ks+1}U^{\ot i}\ot\Ker\pi\ot U^{\ot(ks+1-i)}$$
by Lemma \ref{right}. Now it is not difficult to show that
$$
\beta_{s+1}(1_U^{\ot(s+1)}\ot \rho_k)\left(\sum\limits_{i=1}^{ks+1}U^{\ot i}\ot\Ker\pi\ot U^{\ot(ks+1-i)}\right)\subset M_{\QQ}\ot_{A_{\QQ}}\Ker\vp_{\QQ}
$$
and $\beta_{s+1}(1_U^{\ot(s+1)}\ot \rho_k)\left(\Ker\pi\ot U^{\ot(ks+1)}\right)\subset\Ker\vp_{\QQ}\ot_{A_{\QQ}} M_{\QQ}$. Thus, Condition 3 holds. Analogously, using the surjectivity of the map $\beta_{2s}$, one can prove Condition 4.
\end{Proof}

 We define $\mathfrak{Trip}(U,W,\pi)=(A_{\QQ},M_{\QQ},\vp_{\QQ})$.
If $(f,g):\QQ=(U,W,\pi)\rightarrow (U',W',\pi')=\QQ'$ is a morphism of quivers with corelations, then it induces the morphism $T(g):T(W)\rightarrow T(W')$ such that $T(g)|_{W^{\ot i}}=g^{\ot i}$ for all $i\ge 0$. It is easy to see that $T(g)(I_{\pi})\subset I_{\pi'}$, i.e. $T(g)$ induces a well defined morphism of $\OO$-algebras $T(g):A_{\QQ}\rightarrow A_{\QQ'}.$ Now, it is easy to see that the map $f:(M_{\QQ})_0=U\rightarrow U'=(M_{\QQ'})_0$ can be proceeded to an $A_{\QQ}$-bimodule homomorphism $T(f):M_{\QQ}\rightarrow M_{\QQ'}$ in a unique way.
Thus, we can define $\mathfrak{Trip}(f,g)=\big(T(g),T(f)\big)$. One can check that $\mathfrak{Trip}(f,g)$ is really a morphism of $s$-homogeneous triples.
The proof of the next proposition is standard.

\begin{prop}\label{RT} $(\mathfrak{TRep},\mathfrak{Trip})$ is a pair of quasi inverse equivalences.
\end{prop}

As in the case of algebras, we have a duality $(-)^!=\mathfrak{Trip}\circ (-)^!\circ\mathfrak{TRep}:\HT(\OO, s)\rightarrow \HT(\OO, s)$. Given $\TT\in\HA(\OO,s)$, we will call the $s$-homogeneous triple $\TT^!$ the {\it $s$-dual triple} for $\TT$.




\section{$s$-homogeneous algebras and $s$-homogeneous triples}

In this section we study the relations between $s$-homogeneous algebras and $s$-homogeneous triples.
Our main theorem follows easily from the discussion of the previous sections.

\begin{theorem}\label{equiv}
The categories $\HA(s,\OO)$ and $\HT(s,\OO)$ are equivalent.
\end{theorem}
\begin{Proof} The assertion of the theorem follows from Lemma \ref{DandK} and Propositions \ref{RA} and \ref{RT}. Really,
\begin{multline*}
\mathfrak{Trip}\circ \mathfrak{Coker}\circ \mathfrak{ARep}:\HA(s,\OO)\rightarrow\HT(s,\OO)\mbox{ and }\\
\mathfrak{Alg}\circ \mathfrak{Ker}\circ \mathfrak{TRep}:\HT(s,\OO)\rightarrow\HA(s,\OO)
\end{multline*}
 are quasi inverse equivalences.
\end{Proof}

Now, it is not difficult to describe $\FF=\mathfrak{Trip}\circ \mathfrak{Coker}\circ \mathfrak{ARep}$ and $\GG=\mathfrak{Alg}\circ \mathfrak{Ker}\circ \mathfrak{TRep}$. We describe them using the notation introduced in Section \ref{sAlg}.

If $\Lambda$ is an object of $\HA(s,\OO)$, then $\FF(\Lambda)=(\Lambda^{(s)},\Lambda^{(s,1)},\phi^{(s)}_{\Lambda})$. If $\Gamma\in\HA(s,\OO)$ and $\theta:\Lambda\rightarrow \Gamma$ is a morphism of graded $\OO$-algebras, then
$\FF(\theta):=(\theta|_{\Lambda^{(s)}},\theta|_{\Lambda^{(s,1)}}):\FF(\Lambda)\rightarrow\FF(\Gamma).$

Given  $(A,M,\vp)\in\HT(\OO, s)$, one has $\GG(A,M,\vp)=T(M_0)/\la (\Ker\vp)_0\ra.$
If $(B,L,\psi)$ is another $s$-homogeneous triple and $(f,g):(A,M,\vp)\rightarrow (B,L,\psi)$ is a morphism of $s$-homogeneous triples, then $\GG(f,g):\GG(A,M,\vp)\rightarrow\GG(B,L,\psi)$ is the map induced by $g|_{M_0}:M_0\rightarrow L_0$.

Now we have several corollaries of Theorem \ref{equiv} and the description of $\FF$ and $\GG$.

\begin{coro}\label{AT}
If $\Lambda\in\HA(s,\OO)$, then $(\Lambda^{(s)},\Lambda^{(s,1)},\phi^{(s)}_{\Lambda})\in\HT(s,\OO)$.
\end{coro}

\begin{coro}\label{main_strip}
If $(A,M,\vp)$ is an $s$-homogeneous triple, then there exists a graded $A$-bimodule $S$ concentrated in degree $0$ and an isomorphism of graded $A$-bimodules $\theta:M^{\ot_A^s}\cong S\oplus A_{>0}(1)$ such that $\vp$ equals to the composition
$$
M^{\ot_A^s}\xrightarrow{\theta} S\oplus A_{>0}(1)\twoheadrightarrow A_{>0}(1)\hookrightarrow A(1),
$$
where the second map is the canonical projection on the second summand and the third arrow is the canonical inclusion.
\end{coro}
\begin{Proof} There exists $\Lambda\in\HA(\OO,s)$ such that $(A,M,\vp)\cong\FF(\Lambda)=(\Lambda^{(s)},\Lambda^{(s,1)},\phi^{(s)}_{\Lambda})$. In the notation of Proposition \ref{shom}, we have $\phi^{(s)}_{\Lambda}=\iota_{s,\Lambda}\phi_s^{s-1,1}(\phi_s^{s-2,1}\ot 1_{\Lambda_{s,1}})\dots(\phi_s^{1,1}\ot 1_{\Lambda_{s,1}}^{\ot(s-2)})$, where $\iota_{s,\Lambda}:A_{>0}\cong\Lambda^{(s,s)}\hookrightarrow \Lambda^{(s)}\cong A$ is the canonical inclusion. Now the assertion of the corollary follows from the second item of Proposition \ref{shom}.
\end{Proof}

\begin{coro}\label{uniq}
Let $(A,M,\vp)$ and $(A,M,\vp')$ be $s$-homogeneous triples. If the $A$-bimodule $A_{>0}(1)$ does not contain nonzero direct summands concentrated in degree $0$, then $(A,M,\vp)\cong(A,M,\vp')$.
\end{coro}
\begin{Proof} By Corollary \ref{main_strip}, there are graded $A$-bimodules $S$ and $S'$ concentrated in degree $0$ and isomorphisms of graded $A$-bimodules $\theta:M^{\ot_A^s}\rightarrow S\oplus A_{>0}(1)$ and $\theta':M^{\ot_A^s}\rightarrow S'\oplus A_{>0}(1)$ such that
$\Ker\vp=\theta^{-1}(S)$ and $\Ker\vp'=(\theta')^{-1}(S')$. Suppose that $(\theta')^{-1}(S')\not\subset\theta^{-1}(S)$. Then $\theta(\theta')^{-1}(S')+S$ is a graded $A$-subbimodule of $S\oplus A_{>0}(1)$ that does not lie in $S$, i.e. $T=(\theta(\theta')^{-1}(S')+S)\cap A_{>0}(1)\not=\{0\}$. Now it is clear that the monomorphism $T\hookrightarrow A_{>0}(1)$ splits and, hence, $A_{>0}(1)$ contains nonzero direct summand concentrated in degree $0$. The obtained contradiction proves that $\Ker\vp'\subset\Ker\vp$. The inverse inclusion can be proved in the same way, i.e. $\Ker\vp'=\Ker\vp$. Thus, $\GG(A,M,\vp)=\GG(A,M,\vp')$. Since $\GG$ is an equivalence, we have $(A,M,\vp)\cong(A,M,\vp')$.
\end{Proof}

The next corollary follows directly from Corollary \ref{uniq}.

\begin{coro}\label{uniqalg}
Suppose that $\Lambda,\Gamma\in\HA(\OO,s)$ are such that $\Lambda^{(s,1)}\cong \Gamma^{(s,1)}$ as $\Lambda^{(s)}$-bimodules, where $\Lambda^{(s)}$--bimodule structure on $\Gamma^{(s,1)}$ is induced by some isomorphism $\Lambda^{(s)}\cong \Gamma^{(s)}$. If $\Lambda^{(s)}_{>0}$ does not contain direct $\Lambda^{(s)}$-bimodule summand concentrated in degree $0$, then $\Lambda\cong\Gamma$.
\end{coro}

\begin{Ex} Suppose that $\Lambda\in\HA(\OO,s)$ and $\Lambda^{(s)}\cong \bS(W)$ is the symmetric algebra of the space $W$. Then, due to Corollary \ref{uniqalg}, $\Lambda$ can be uniquely recovered from the $\bS(W)$-bimodule $\Lambda^{(s,1)}$.
\end{Ex}

\begin{Ex} Suppose that $\Lambda\in\HA(\OO,s)$ and $\Lambda^{(s)}\cong \bL(W)$ is the exterior algebra of the space $W$ with $\dim_{\kk}W\ge 2$. Then, due to Corollary \ref{uniqalg}, $\Lambda$ can be uniquely recovered from the $\bL(W)$-bimodule $\Lambda^{(s,1)}$. On the other hand, we will show later that it is not true in the case $\dim_{\kk}W=1$.
\end{Ex}

Note that the functors $\FF$ and $\GG$ respect the duality $(-)^!$, i.e., for $\Lambda\in\HA(\OO, s)$ and $\TT\in\HT(\OO, s)$, one has
$\FF(\Lambda^!)\cong \FF(\Lambda)^!$ and $\GG(\TT^!)\cong \GG(\TT)^!$.

\section{$s$-Koszulity}

In this section we discuss the notion of an $s$-Koszul algebra.

\begin{Def} {\rm The $s$-homogeneous algebra $\Lambda$ is called {\it $s$-Koszul} if $\Ext_{\Lambda}^i(\OO,\OO)$ is concentrated in degree $-\chi_s(i)$, where
\begin{equation}
\chi_s(i)=\begin{cases}
\frac{is}{2},&\mbox{if $2\mid i$},\\
\frac{(i-1)s}{2}+1,&\mbox{if $2\nmid i$}.
\end{cases}
\end{equation}
The $2$-Koszul algebras are called simply {\it Koszul algebras}.
}
\end{Def}

Let $\Lambda$ be an $s$-homogeneous algebra and $(V,U,\iota)=\mathfrak{ARep}(\Lambda)$. Let us define the components of the graded vector space $R$ by the equality $R_n=\cap_{i+s+j=n}U^{\ot i}\ot \Im\iota\ot U^{\ot j}$.
The inclusion $R_{n+m}\hookrightarrow R_n\ot R_m$ induces a map $R_n^*\ot R_m^*\rightarrow R_{n+m}^*$ that gives a graded $\OO$-algebra structure on $R^*=\oplus_{n\ge 0}R_n^*$. In fact, it is well known that $R^*\cong \Lambda^!$. We also introduce the graded space $I$ with $n$-th component $I_n=\sum_{i+s+j=n}U^{\ot i}\ot \Im\iota\ot U^{\ot j}$. Note that, in fact, $I$ is an ideal in $T(U)$ such that $\GG(V,U,\iota)=T(U)/I$.

For $n>m$, the inclusion $R_n\hookrightarrow R_m\otimes U^{\ot(n-m)}$ induces a $\Lambda$-module homomorphism $d_m^n:R_n\otimes\Lambda\rightarrow R_m\otimes\Lambda$. Then we can define the {\it generalized Koszul complex} $K$ of $\Lambda$ in the following way. Its $n$-th member is $K_n=R_{\chi_s(n)}\otimes\Lambda$. The differential is defined by the equality $d(K)_n=d^{\chi_s(n+1)}_{\chi_s(n)}:K_{n+1}\rightarrow K_n$. Note that there is a surjective homomorphism $\mu_K:K_0\rightarrow \OO$ induced by the composition of the isomorphism $\OO\otimes\Lambda\cong\Lambda$ and the canonical projection $\Lambda\twoheadrightarrow\Lambda_0=\OO$. It is proved in \cite{B} that $\Lambda$ is $s$-Koszul if and only if $K$ is exact in positive degrees (though $\OO=\kk$ there, it is not difficult to transfer the arguments of Berger to our case). Direct calculations (which are fulfilled in \cite{B} for $\OO=\kk$) show that $K$ is exact in positive degrees if and only if the following conditions are satisfied for all $k,n\ge 0$:
\begin{equation}\label{main_cnd}
\begin{array}{c}
R_{ns+1}\otimes U^{\ot k}\cap U^{\ot(ns)}\ot I_{k+1}=R_{(n+1)s}\ot U^{\ot(k-s+1)}+R_{ns+1}\ot I_k,\\
R_{(n+1)s}\otimes U^{\ot k}\cap U^{\ot(ns+1)}\ot I_{k+s-1}=R_{(n+1)s+1}\ot U^{\ot(k-1)}+R_{(n+1)s}\ot I_k.
\end{array}
\end{equation}

It is shown in \cite{B} that \eqref{main_cnd} is satisfied if and only if
\begin{equation}\label{ec}
I_s\ot U^{\ot(s-1)}\cap U\ot I_{2s-2}=R_{s+1}\ot U^{\ot(s-2)}
\end{equation}
and  for all $k,n\ge 0$ we have
\begin{equation}\label{dc}
\begin{array}{c}
R_{ns+1}\otimes U^{\ot k}\cap U^{\ot(ns)}\ot I_{k+1}=R_{ns+1}\ot U^{\ot k}\cap U^{\ot (ns)}\ot I_s\ot U^{\ot(k-s+1)}+R_{ns+1}\ot I_k,\\
\begin{aligned}
R_{(n+1)s}\otimes U^{\ot k}&\cap U^{\ot(ns+1)}\ot I_{k+s-1}\\
&=R_{(n+1)s}\ot U^{\ot k}\cap U^{\ot(ns+1)}\ot I_{2s-2}\ot U^{\ot(k-s+1)}+R_{(n+1)s}\ot I_k.
\end{aligned}
\end{array}
\end{equation}
We call the condition \eqref{ec} the {\it extra condition} and call the conditions \eqref{dc} the {\it distributivity conditions}.

If $s=2$, then it is well known that $\Lambda$ is Koszul if and only if $\Lambda^!$ is. For $s>2$ there are examples where it is not so. One can show that the conditions for the $s$-Koszulity of the algebra $\Lambda^!$ are conditions \eqref{ec} and \eqref{dc} with $R$ and $I$ interchanged.

If $s=2$ and $\Lambda$ is Koszul, then $\Ext_{\Lambda}^*(\OO,\OO)\cong (\Lambda^!)^{\rm op}$ as a graded algebra. For $s>2$ the situation changes a little.
Let $\FF(\Lambda^!)=(A,M,\vp)$. It is proved in \cite{GMMZ} that if $\Lambda$ is $s$-Koszul, then $A$ is a Koszul algebra and $M$ is linear as left and right $A$-module.
If $\Lambda$ is $s$-Koszul and $s>2$, then $\Ext_{\Lambda}^*(\OO,\OO)\cong A\ltimes M$ as an algebra, where $A\ltimes M$ is the trivial extension of $A$ by $M$, i.e. its underlying space is $A\oplus M$ and the multiplication is given by the equality $(a,x)(b,y)=(ab,ay+xb)$ for $a,b\in A$ and $x,y\in M$. If we define the grading on $A\ltimes M$ by the equalities $(A\ltimes M)_{2n}=A_n$ and $(A\ltimes M)_{2n+1}=M_n$ for $n\ge 0$, then the isomorphism above will become degree preserving.
From our previous results we get the following easy corollary.

\begin{coro}
Suppose that $\Lambda$ and $\Gamma$ are $s$-Koszul algebras such that $\Ext^*_{\Lambda}(\OO,\OO)\cong\Ext^*_{\Gamma}(\OO,\OO)$. If $\Lambda\not\cong\Gamma$, then there exists nonzero $\theta\in\Ext^2_{\Lambda}(\OO,\OO)$ such that $\theta\Ext^2_{\Lambda}(\OO,\OO)=\Ext^2_{\Lambda}(\OO,\OO)\theta=0$.
\end{coro}
\begin{Proof} Let us set $(A,M,\vp)=F(\Lambda^!)$ and $(B,L,\psi)=F(\Gamma^!)$. Then $\Ext^*_{\Lambda}(\OO,\OO)\cong A\ltimes M$ and $\Ext^*_{\Gamma}(\OO,\OO)\cong B\ltimes L$ with the gradings defined above. By the degree argument this means that there is an isomorphism of graded algebras $A\cong B$ and an isomorphism of graded $A$-bimodules $M\cong L$, where the $A$-bimodule structure on $L$ is induced by the just mentioned isomorphism of algebras. Then $(B,L,\psi)\cong (A,M,\vp')$ for some $\vp':M^{\ot_A^s}\rightarrow A(1)$. If $\Lambda\not\cong\Gamma$, then we have $(A,M,\vp)\not\cong (A,M,\vp')$, and hence $A_{>0}(1)$ contains nonzero $A$-bimodule summand concentrated in degree zero by Corollary \ref{uniq}, i.e. there is a nonzero element $a\in A_1$ such that $aA_1=A_1a=0$. Since $A_1$ can be identified with $\Ext^2_{\Lambda}(\OO,\OO)$, this is exactly the required assertion.
\end{Proof}

\section{Hilbert series}

Let us now discuss the notion of Hilbert series and its relations with the notion of $s$-Koszulity. For an $\OO$-bimodule $W$ we will denote by $\dim_W$ the endomorphism of $\OO$ defined by the equality $\dim_W(e_j)=\sum\limits_{i=1}^D\big(\dim_{\kk}e_iWe_j\big)e_i$ for $1\le j\le D$.

\begin{Def} {\rm Let $W=\oplus_{k\ge 0}W_k$ be a nonnegatively graded $\OO$-bimodule. The {\it Hilbert series} of $W$ is the map $\HH_W(t):\OO\rightarrow \OO[[t]]$ defined by the equality
$$\HH_W(t)=\sum\limits_{k=0}^{\infty}t^k\dim_{W_k}.$$
}
\end{Def}

If $\Lambda\in\HA(\OO,s)$ is $s$-Koszul, then, by the results of \cite{BS}, one has
$$
\Big(\HH_A(t^s)-t\HH_M(t^s)\Big)\HH_{\Lambda}(t)=\Id_{\OO},
$$
where $(A,M,\vp)=\FF(\Lambda^!)$. To prove this it is enough to note that the left part of the equality is $\sum\limits_{k=0}^{\infty}(-1)^k\HH_{K_k}(t)$, where $K$ is the generalized Koszul complex of $\Lambda$. We will denote by $O(t^n)$ the set $t^n\Hom_{\kk}(\OO, \OO[[t]])\subset\Hom_{\kk}(\OO, \OO[[t]])$.
In what follows, $f=g+O(t^n)$ means $f-g\in O(t^n)$. Then we have the following result

\begin{lemma}\label{ECVH} $\Lambda\in\HA(\OO,s)$ satisfies the extra condition if and only if $$\Big(\HH_A(t^s)-t\HH_M(t^s)\Big)\HH_{\Lambda}(t)=\Id_{\OO}+O(t^{2s}).$$
\end{lemma}
\begin{Proof} Let us set $(V,U,\iota)=\mathfrak{ARep}(\Lambda)$,
$$I_n=\sum_{i+s+j=n}U^{\ot i}\ot \Im\iota\ot U^{\ot j}\mbox{, and }R_n=\cap_{i+s+j=n}U^{\ot i}\ot \Im\iota\ot U^{\ot j}.$$
Note that $\HH_A(t^s)=\Id_{\OO}+t^s\dim_{I_s}+O(t^{2s})$, $t\HH_M(t^s)=t\dim_U+t^{s+1}\dim_{R_{s+1}}+O(t^{2s})$, and
$$
\HH_{\Lambda}(t)=\sum\limits_{n=0}^{s-1}t^n\dim_{U^{\ot n}}+\sum\limits_{n=s}^{2s-1}t^n(\dim_{U^{\ot n}}-\dim_{I_n})+O(t^{2s}).
$$
Thus, we have
\begin{multline*}
\Big(\HH_A(t^s)-t\HH_M(t^s)\Big)\HH_{\Lambda}(t)\\
=(\Id_{\OO}-t\dim_U+t^s\dim_{I_s}-t^{s+1}\dim_{R_{s+1}})\left(\sum\limits_{n=0}^{2s-1}t^n\dim_{U^{\ot n}}-\sum\limits_{n=s}^{2s-1}t^n\dim_{I_n}\right)+O(t^{2s})\\
=\Id_{\OO}+\sum\limits_{n=s}^{2s-1}t^n(\dim_{I_s\ot U^{\ot (n-s)}}-\dim_{I_n})+\sum\limits_{n=s+1}^{2s-1}t^n(\dim_{U\ot I_{n-1}}-\dim_{R_{s+1}\ot U^{\ot (n-s-1)}})+O(t^{2s})\\
=\Id_{\OO}+\sum\limits_{n=s+1}^{2s-1}t^n(\dim_{I_s\ot U^{\ot (n-s)}\cap U\ot I_{n-1}}-\dim_{R_{s+1}\ot U^{\ot (n-s-1)}})+O(t^{2s}).
\end{multline*}
Since $R_{s+1}\ot U^{\ot (n-s-1)}\subset I_s\ot U^{\ot (n-s)}\cap U\ot I_{n-1}$, the equality from the assertion of the lemma is satisfied if and only if
$R_{s+1}\ot U^{\ot (n-s-1)}= I_s\ot U^{\ot (n-s)}\cap U\ot I_{n-1}$ for any $s+1\le n\le 2s-1$. It is clear that the last mentioned equality holds for $s+1\le n\le 2s-1$ if and only if it holds for $n=2s-1$, i.e. if and only if the condition \eqref{ec} is satisfied.
\end{Proof}

Note now that the Hilbert series of $\Lambda$ are fully determined by the first two components of $\FF(\Lambda)$.

\begin{lemma}\label{HSWT} Let $\Lambda\in\HA(\OO,s)$ and $\FF(\Lambda)=(B,L,\psi)$. Then $\HH_{\Lambda}(t)=\sum\limits_{k=0}^{s-1}t^k\HH_{L^{\ot_B^k}}(t^s)$.
\end{lemma}
\begin{Proof} Follows from the fact that the map $(\Lambda^{(s,1)})^{\ot_{\Lambda^{(s)}}^k}\rightarrow \Lambda^{(s,k)}$ is bijective for $0\le k\le s-1$ by the second item of Proposition \ref{shom}.
\end{Proof}

\begin{coro}\label{ecc} Let $\Lambda\in\HA(\OO,s)$ and $\FF(\Lambda)=(B,L,\psi)$. Then $\Lambda$ satisfies the extra condition if and only if
\begin{multline*}
\Big(\Id_{\OO}-t\dim_{L_0}+t^s(\dim_{L_0^{\ot s}}-\dim_{B_1})\\
-t^{s+1}(\dim_{L_0^{\ot(s+1)}\oplus L_1}-\dim_{L_0\ot B_1\oplus B_1\ot L_0})\Big)\sum\limits_{k=0}^{s-1}t^k\HH_{L^{\ot_B^k}}(t^s)-\Id_{\OO}\in O(t^{2s}),
\end{multline*}
In particular, if the algebra $\Lambda^{(s)}$ and the $\Lambda^{(s)}$-bimodule $\Lambda^{(s,1)}$ are known, then it is known if $\Lambda$ satisfies the extra condition or not.
\end{coro}
\begin{Proof} We will be free to use the notation of the proof of Lemma \ref{ECVH}. Since
$$\dim_{I_s}=\dim_{\Lambda_1^{\ot s}}-\dim_{\Lambda_s}=\dim_{L_0^{\ot s}}-\dim_{B_1},$$ we have $\HH_A(t^s)=\Id_{\OO}+t^s(\dim_{L_0^{\ot s}}-\dim_{B_1})+O(t^{2s}).$

The exact sequence $R_{s+1}\hookrightarrow U\ot I_s\oplus I_s\ot U \twoheadrightarrow I_{s+1}$ gives the equality
\begin{multline*}
\dim_{R_{s+1}}=\dim_{U\ot I_s\oplus I_s\ot U}-\dim_{I_{s+1}}\\
=2\dim_{L_0^{\ot(s+1)}}-\dim_{L_0\ot B_1\oplus B_1\ot L_0}-(\dim_{L_0^{\ot(s+1)}}-\dim_{L_1})
=\dim_{L_0^{\ot(s+1)}\oplus L_1}-\dim_{L_0\ot B_1\oplus B_1\ot L_0}.
\end{multline*}
Since $L_0=U$, we have
$t\HH_M(t^s)=t\dim_{L_0}+t^{s+1}(\dim_{L_0^{\ot(s+1)}\oplus L_1}-\dim_{L_0\ot B_1\oplus B_1\ot L_0})+O(t^{2s}).$

Now the assertion of the corollary follows from Lemmas \ref{ECVH} and \ref{HSWT}.
\end{Proof}

\section{Examples and applications}

In this section we give some examples showing how the technique of $s$-homogeneous triples works. In particular, we will discuss some of the results of \cite{B2}.

In this section we set $\OO=\kk$. Moreover, we assume for simplicity that $\kk$ is algebraically closed.
Let, as before, $\Lambda$ be an $s$-homogeneous algebra and $(A,M,\vp)=\FF(\Lambda^!)$.

\subparagraph{Example 1, $A=\kk[x]$.} As it was mentioned before, in this case $\Lambda^!$ and $\Lambda$ are determined by the bimodule $M$. By Corollary \ref{main_strip}, we have $M^{\ot_A^s}\cong \kk^l\oplus A_{>0}=\kk^l\oplus \kk[x]$ for some $l\ge 0$. From this condition and the fact that $M$ has to have linear presentation, it is not difficult to deduce that $M\cong \kk^m\oplus \kk[x]_{\alpha}$, where $m\ge 0$ and $\alpha$ is an automorphism of $\kk[x]$ sending $x$ to $\epsilon x$ for some $s$-th root of unit $\epsilon$. Thus, $l=(m+1)^s-1$ and $\vp$ is the canonical projection
$M^{\ot_A^s}\cong\kk^l\oplus (\kk[x]_{\alpha})^{\ot_{\kk[x]}^s}\twoheadrightarrow (\kk[x]_{\alpha})^{\ot_{\kk[x]}^s}\cong\kk[x]$. It is easy to verify that the third condition from the definition of an $s$-homogeneous triple holds if and only if $\epsilon=1$. In this case $\Lambda=\kk\langle x_1,\dots,x_m,y\rangle/(y^s)$ is $s$-Koszul and
it is not difficult to see that $\Lambda^!$ is $s$-Koszul if and only if either $s=2$ or $m=0$.

\subparagraph{Example 2, $A=\kk[x]/(x^2)$.} It is not difficult to show that in this case $M\cong \kk^m$ for some $m\ge 2$. Then $(\Lambda^!)_n=0$ for $n>s$ and $\Lambda$ is $s$-Koszul if and only if $\Lambda$ has global dimension $2$, i.e. the generalized Koszul complex of $\Lambda$ is exact in the second term.
It follows from the results of \cite{B2} that $\Lambda$ is $s$-Koszul if and only if it satisfies the extra condition, i.e. if and only if $(1-mt+t^s)\HH_{\Lambda}(t)-1\in O(t^{2s})$. It is easy to see that the situation under consideration occurs if and only if $\Lambda=\kk\langle x_1,\dots,x_m\rangle/(f)$, where $f$ is some homogeneous polynomial in $x_1,\dots,x_m$ of degree $s$ such that $f\not=g^s$ for any linear polynomial $g$. For example, $\Lambda$ is $s$-Koszul for $f=x_1x_2^{s-1}$ and $\Lambda$ is not $s$-Koszul for $f=x_1x_2^{s-2}x_1$. More detailed description of the situation is given in the next proposition.

\begin{prop}
Suppose that $\Lambda=\kk\langle x_1,\dots,x_m\rangle/(f)$, where $f$ is some homogeneous polynomial in $x_1,\dots,x_m$ of degree $s$. Then $\Lambda$ is $s$-Koszul if and only if one of the following two conditions holds:
\begin{enumerate}
\item  $f=g^s$ for some linear polynomial $g$;
\item if $f=gh_1=h_2g$ for some polynomials $g$, $h_1$ and $h_2$, then $\deg g\in\{0,\deg f\}$.
\end{enumerate}
\end{prop}
\begin{Proof} If the first condition holds, then the $s$-Koszulity follows from the argument above.

If the first condition does not hold, then it is easy to see that the generalized Koszul complex of $\Lambda$ has the form
$$
\kk f\ot\Lambda\xrightarrow{d_1}\oplus_{i=1}^m\kk x_i\ot\Lambda\xrightarrow{d_0}\Lambda(\twoheadrightarrow \kk),
$$
where $d_0(x_i\ot 1)=x_i$ ($1\le i\le m$) and $d_1(f\ot 1)=(x_1\ot f_1,\dots,x_m\ot f_m)$ for such $f_1,\dots, f_m$ that $f=\sum\limits_{i=1}^mx_if_i$. Then $\Lambda$ is $s$-Koszul if and only if $d_1$ is injective.
If $gh_1=h_2g$ for some $g$, $h_1$ and $h_2$ such that $0<\deg g<\deg f$, then it is easy to see that $f\ot h_1\not=0$ and $d_1(f\ot h_1)=0$, i.e. $d_1$ is not injective.

Suppose that the second condition holds and $d_1(f\ot h)=0$ for some polynomial $h$. This means that $fh=\sum\limits_{j=1}^ku_{1,j}fu_{2,j}$ in $\kk\langle x_1,\dots,x_m\rangle$ for some  $u_{1,j},u_{2,j}\in \kk\langle x_1,\dots,x_m\rangle$ such that $\deg u_{1,j}>0$. Let us consider the map $T:\kk\langle x_1,\dots,x_m\rangle\rightarrow \kk\langle x_1,\dots,x_m\rangle$ defined by the equality $T(y)=\sum\limits_{j=1}^ku_{1,j}yu_{2,j}-yh$. Since $T(f)=0$ we have $$T(y)=\sum\limits_{j=1}^r(v_{1,j}fv_{2,j}yv_{3,j}-v_{1,j}yv_{2,j}fv_{3,j})$$ for some $v_{1,j},v_{2,j},v_{3,j}\in \kk\langle x_1,\dots,x_m\rangle$ by \cite[Corollary 1.6]{Gera}. Now it is easy to see that $h$ has the form $h=\sum\limits_{j=1}^mh_{1,j}fh_{2,j}$ for some $h_{1,j},h_{2,j}\in \kk\langle x_1,\dots,x_m\rangle$. Thus, $f\ot h=0$ in $\kk f\ot\Lambda$, and hence $d_1$ is injective.
\end{Proof}

On the other hand, we have
\begin{multline*}
\Big(1-t\dim_{M_0}+t^s(\dim_{M_0^{\ot s}}-\dim_{A_1})\\
-t^{s+1}(\dim_{M_0^{\ot(s+1)}\oplus M_1}-\dim_{M_0\ot A_1\oplus A_1\ot M_0})\Big)\sum\limits_{k=0}^{s-1}t^k\HH_{M^{\ot_A^k}}(t^s)\\
=\big(1-mt+(m^s-1)t^s-(m^{s+1}-2m)t^{s+1}\big)\left(\sum\limits_{i=0}^{s-1}(mt)^i+t^s\right)=1+t^s\sum\limits_{i=2}^{s-1}(mt)^i+O(t^{2s}),
\end{multline*}
i.e. $\Lambda^!$ does not satisfy the extra condition by Corollary \ref{ecc} if $s>2$.

\begin{coro} Let $\Lambda\in\HA(\kk,s)$ ($s\ge 3$) be such that $\dim_{\kk}\Lambda_s=1$. Then $\Lambda$ is $s$-Koszul if and only if
$\Lambda\cong\kk[x]$ .
\end{coro}

Now we  give an example where the  algebra $A$ has two generators. As a first step in this direction, we will get a restriction on the number of relations for such an algebra.

\begin{theorem} Let $(A,M,\vp)$ be an $s$-homogeneous triple over $\kk$. If $A=\kk\langle x,y\rangle/I$ for some quadratic ideal $I$, then $\dim_{\kk}A_2\le 2$.
\end{theorem}
\begin{Proof} Let $(U,W,\pi)=\mathfrak{TRep}(A,M,\vp)$, i.e. $U=M_0$, $W=A_1$, and $\pi=\vp_0$. By our assumption, we have $\dim_{\kk}\Im\pi=\dim_{\kk}W=2$. Suppose that $\dim_{\kk}A_2> 2$. This means that
$$\dim_{\kk}(\pi\ot \pi)\left(\sum\limits_{i=0}^sU^{\ot i}\ot\Ker\pi\ot U^{\ot(s-i)}\right)\le 1.$$
Since $\Ker(\pi\ot\pi)=\Ker\pi\ot U^{\ot s}+U^{\ot s}\ot\Ker\pi$, we have
$$\sum\limits_{i=0}^sU^{\ot i}\ot\Ker\pi\ot U^{\ot(s-i)}\subset \Ker\pi\ot U^{\ot s}+U^{\ot s}\ot\Ker\pi+V$$
for some $V\subset U^{\ot(2s)}$ of dimension not more than $1$. For $S\subset U^{\ot l}$ we as usually introduce $S^{\perp}=\{\alpha\in (U^*)^{\ot l}\mid \alpha(u)=0\,\,\forall u\in S\}$. Then we have
\begin{multline*}
(\Ker\pi)^{\perp}\ot (U^*)^{\ot s}\cap(U^*)^{\ot s}\ot(\Ker\pi)^{\perp}\cap V^{\perp}=(\Ker\pi\ot U^{\ot s}+U^{\ot s}\ot\Ker\pi+V)^{\perp}\\
\subset \left(\sum\limits_{i=0}^sU^{\ot i}\ot\Ker\pi\ot U^{\ot(s-i)}\right)^{\perp}=\bigcap\limits_{i=0}^s(U^*)^{\ot i}\ot(\Ker\pi)^{\perp}\ot (U^*)^{\ot(s-i)}.
\end{multline*}
Let us fix some basis $x_1,\dots,x_n$ of $U^*$. Note that any element of $(U^*)^{\ot l}$ can be written as a linear combination of words of length $l$ in letters $x_1,\dots,x_n$, i.e. as an element of $T_n$, where $T=\kk\langle x_1,\dots,x_n\rangle$. Since $\dim_{\kk}(\Ker\pi)^{\perp}=2$, there are two element $f_1,f_2\in T_s$ that form a basis of $(\Ker\pi)^{\perp}$. The condition above can be rewritten in the form
$$
(\kk f_1+\kk f_2)T_s\cap T_s(\kk f_1+\kk f_2)\cap L\subset \bigcap\limits_{i=0}^sT_i(\kk f_1+\kk f_2)T_{s-i},
$$
where $L\subset T_{2s}$ is some subspace of dimension not less than $n^{2s}-1$. Let us introduce the order $x_1>x_2>\cdots>x_n$ on variables and the lexicographic order on the set of monomials of the same length, i.e. $w_1>w_2$ for monomials $w_1,w_2$ of the same length if and only if there are monomials $w, w_1',w_2'$ and integers $1\le i<j\le n$ such that $w_1=wx_iw_1'$ and $w_2=wx_jw_2'$. As usually, for $f\in T_l$ we denote by $tip(f)$ the biggest monomial that has a nonzero coefficient in the decomposition of $f$. Clearly, we may assume that $tip(f_1)\not=tip(f_2)$. Then it is easy to see that
$tip(w)\in T_i\big(\kk \,tip(f_1)+\kk \,tip(f_2)\big)T_{s-i}$ for any $w\in T_i(\kk f_1+\kk f_2)T_{s-i}$.

Since $\dim_{\kk}L\ge n^{2s}-1$ and the monomials $tip(f_1)tip(f_1)$, $tip(f_1)tip(f_2)$, $tip(f_2)tip(f_1)$, and $tip(f_2)tip(f_2)$ are pairwise not equal, it is easy to see that at least three of the four listed monomials belong to $T_i\big(\kk \,tip(f_1)+\kk \,tip(f_2)\big)T_{s-i}$ for any $0\le i\le s$.
It is clear from our argument that we may assume that $f_1$ and $f_2$ are monomials.

Since one of the monomials $f_1f_1$ and $f_1f_2$ belongs to $T_i\big(\kk \,tip(f_1)+\kk \,tip(f_2)\big)T_{s-i}$ for any $0\le i\le s$, it is clear that, for some letter $x$, there are maximum two different words of the form $x_ix_j$ that occur in $f_1x$. Analogous assertion can be proved for some word of the form $yf_1$, where  $y$ is a letter. Then we may assume that either $f_1=x_1^s$ or $f_1=(x_1x_2)^tx_1^{s-2t}$, where $t=\upper{\frac{s-1}{2}}$. Analogous arguments show that we may set $f_2=x_2^s$ in the first case and $f_2=(x_2x_1)^tx_2^{s-2t}$ in the second case. Direct verifications show that all the obtained pairs $(f_1,f_2)$ do not satisfy the required conditions.
\end{Proof}

Thus, it makes sense to consider $A=\kk\langle x,y\rangle/I$, where $I$ is generated by two, three or four quadratic relations.

It is clear that the only case of four relations is $A=\kk\langle x,y\rangle/(x,y)^2$. Note that $s$-Koszul algebras $\Lambda$ with $\FF(\Lambda^!)$ of the form $(\kk\langle x_1,\dots,x_n\rangle/(x_1,\dots,x_n)^2,M,\vp)$ correspond exactly to local $s$-Koszul algebras with $n$ relations of global dimensions $2$ and $3$.

It is not difficult to show that in the case of three relations $A$ is isomorphic to one of the following algebras:
$$\kk\langle x,y\rangle/(x^2,xy,yx),\,\kk\langle x,y\rangle/(x^2,y^2,xy+qyx)\,\,(q\in\kk),\mbox{ and }\kk\langle x,y\rangle/(x^2,xy+yx,xy+y^2).$$
Note that $\kk\langle x,y\rangle/(x^2,y^2,xy+qyx)\cong \kk\langle x,y\rangle/(x^2,y^2,xy+q'yx)$ for $q\not=q'$ if and only if $qq'=1$.
It is not difficult to show, using the results of \cite{B0}, that in the case of two relations $A$ is Koszul if and only if it is isomorphic to one of the following algebras:
$$
\kk\langle x,y\rangle/(xy,yx),\,\kk\langle x,y\rangle/(x^2,y^2)\,,\kk\langle x,y\rangle/(x^2,yx),\,\mbox{ and }\kk\langle x,y\rangle/(x^2,xy+qyx)\,(q\in\kk).
$$
Note that all the listed algebras except $\kk\langle x,y\rangle/(x^2,xy,yx)$ and $\kk\langle x,y\rangle/(x,y)^2$ satisfy the condition of Corollary \ref{uniq}. Thus, if $\dim_{\kk}A_1=2$, then $\Lambda$ can be recovered from the pair $(A,M)$ except the cases $A=\kk\langle x,y\rangle/(x^2,xy,yx)$ and $A=\kk\langle x,y\rangle/(x,y)^2$.

As an example of an application of our technique, we consider in this paper the algebras $\kk\langle x,y\rangle/(xy,yx)$ and $\kk\langle x,y\rangle/(x^2,y^2)$. We believe that our technique can be applied to other cases to obtain a classification of $s$-Koszul algebras with two homogeneous relations and of $s$-Koszul algebras with two dimensional $s$-th component.

Note that, for any quadratic algebra $A=\kk\langle x,y\rangle/I$, any finitely generated right $A$-module $M$ linear until the first degree is isomorphic to a direct sum of indecomposable modules linear until the first degree. At the same time, any indecomposable $A$-module linear until the first degree is isomorphic to either a module of the form
$$A^n/\langle\{f_ix-f_{i+1}y\}_{1\le i\le n-1},af_1y, bf_nx\rangle_A\,\,\,(n\ge 1, a,b\in\{0,1\})$$
or a module of the form
$$A^n/\langle\{f_ix-f_{i+1}y\}_{1\le i\le n-1},f_nx-q f_1y\rangle_A\,\,\,(n\ge 1,q\in\kk^*).$$
Here and further $f_1,\dots,f_n$ denote standard generators of the free module $A^n$.

\subparagraph{Example 3, $A=\kk\langle x,y\rangle/(xy,yx),\kk\langle x,y\rangle/(x^2,y^2)$.} As usually, for $a\in A$, we will denote by $aA$, $Aa$, and $AaA$ respectively the right submodule, left submodule and subbimodule of $A$ generated by $a$. In all cases we have an isomorphism $A_{>0}\cong xA\oplus yA$ of graded right modules and an isomorphism $A_{>0}\cong Ax\oplus Ay$ of graded left modules.

Suppose $(A,M,\vp)$ is an $s$-homogeneous triple. In particular, $M$ is linear until the first degree as left and right $A$-module and satisfies the condition $M^{\ot_A^s}\cong A_{>0}(1)\oplus S$ for some $A$-bimodule $S$ concentrated in zero degree by Corollary \ref{main_strip}.

We start by proving that $M\cong S_r\oplus xA(1)\oplus yA(1)$ as a graded right $A$-module for some module $S_r$ concentrated in zero degree. Let $J$ denote the ideal $A_{>0}$ of $A$. For a subset $I\in A$ and a right $A$-module $X$ let us introduce $l\Ann_IX=\{u\in X\mid uI=0\}$.
Note that $l\Ann_{J^2}M^{\ot_A^s}=l\Ann_{J}M^{\ot_A^s}=S$. Suppose that $1\le k< s$ and $u\in l\Ann_{J^2}M^{\ot_A^k}\setminus l\Ann_{J}M^{\ot_A^k}$. We have
\begin{multline*}
\left(u\ot_AM^{\ot_A^{s-k}}\right)J^2=\left(u\ot_AM^{\ot_A^{s-k}}\right)\vp(M^{\ot_A^s})\vp(M^{\ot_A^s})\\
=u\vp(M^{\ot_A^s})\vp(M^{\ot_A^s})\ot_AM^{\ot_A^{s-k}}=uJ^2\ot_AM^{\ot_A^{s-k}}=0.
\end{multline*}
Hence, $u\ot_AM^{\ot_A^{s-k}}\subset S$, i.e. $\vp(u\ot_AM^{\ot_A^{s-k}})=0$. On the other hand, we have $0\not=uJ=u\vp(M^{\ot_A^s})=\vp(u\ot_AM^{\ot_A^{s-k}})M^{\ot_A^k}$. The obtained contradiction shows that $l\Ann_{J^2}M^{\ot_A^k}= l\Ann_{J}M^{\ot_A^k}$ for any $1\le k< s$.

According to the classification of indecomposable $A$-module linear until the first degree given above, direct right $A$-module summands of $M^{\ot_A^k}$ ($1\le k\le s$) can have the following forms:
$$
\kk,\,A,\,xA(1),\,yA(1),\,N=(A\oplus A)/\langle(x,-y)\rangle_A.
$$
It is easy to show by induction that if $M$ has a direct right $A$-module summand isomorphic to $A$, then $M^{\ot_A^k}$ has such a summand for any $k\ge 1$. On the other hand, $M^{\ot_A^s}$ does not have direct summand isomorphic to $A$.
Note that $M^{\ot_A^k}$ has  a direct summand isomorphic to $N$ if and only if $(l\Ann_{J}M^{\ot_A^k})_1\not=0$. Suppose that $M$ has a direct summand isomorphic to $N$.
Since $J(l\Ann_{J}M)_1\subset MJ^2$, $\big(J(l\Ann_{J}M)_1\big)J=0$ and $l\Ann_{J}(MJ^2)=0$, we have $J(l\Ann_{J}M)_1=0$. The argument contained in the previous part of the proof shows that $M$ has a direct left $A$-module summand isomorphic to $N'=(A\oplus A)/\langle(x,-y)\rangle_{A^{\rm op}}$.
Direct calculations show that $\Tor_A(xA(1),N')=\Tor_A(yA(1),N')=0$. Applying the functor $-\ot_AM$ to the short exact sequence $\kk(-1)\hookrightarrow N \twoheadrightarrow xA(1)\oplus yA(1)$ we get the long exact sequence
$$
\cdots\rightarrow \Tor(xA(1)\oplus yA(1),M)\xrightarrow{\beta} (M/JM)(-1)\xrightarrow{\alpha} N\ot_A M \twoheadrightarrow \big(xA(1)\oplus yA(1)\big)\ot_A M.
$$
The argument above shows that $(M/JM)(-1)$ has a direct summand isomorphic to $(N'/JN')(-1)$ that does not belong to the image of $\beta$. Thus, $\beta$ is not surjective and $\Im\alpha\subset (l\Ann_{J}N\ot_A M)_1$ is nonzero.
Since $N\ot_A M$ is a direct summand of $M\ot_A M$, the argument above shows that $N\ot_A M$ has a direct summand isomorphic to $N$. Then it is easy to show by induction that $M^{\ot_A^k}$ has a summand isomorphic to $N$ for any $k\ge 1$. The obtained contradiction shows that
$M\cong S_r\oplus \big(xA(1)\big)^{k_1}\oplus \big(yA(1)\big)^{k_2}$ as a graded right $A$-module for some module $S_r$ concentrated in zero degree and integers $k_1,k_2\ge 0$.

Since $Mx\not=0$ and $My\not=0$, it is easy to see that $k_1,k_2\ge 1$. Suppose that $k_1>1$. Then $\left(xA(1)\ot_AM^{\ot_A^{s-1}}\right)^{k_1}$ is a direct summand of $M^{\ot_A^s}$, and hence $xA(1)\ot_AM^{\ot_A^{s-1}}$ is concentrated in zero degree.
Then $k_2=1$ and there exists a subspace $U'\subset M_0$ of codimension $1$ such that $J(U'\ot_AM^{\ot_A^{s-1}})=0$. Applying the dual argument we get a subspace $U\subset M_0$ of codimension $1$ such that $(M^{\ot_A^{s-1}}\ot_AU)J=0$.
Now it follows from the direct sum right $A$-module decomposition of $M$ that either $M^{\ot_A^s}x=0$ or $M^{\ot_A^s}y=0$. The obtained contradiction shows that $M\cong S_r\oplus xA(1)\oplus yA(1)$ as a graded right $A$-module and finishes our first step.

Since $JS_r\subset MJ$, $(JS_r)J=0$, and $l\Ann_J(MJ)=0$, we have  $JS_r=0$, i.e. $S_r=\kk^m$ is a direct $A$-bimodule summand of $M$ concentrated in degree $0$.
It remains to describe left $A$-module structure on $xA(1)\oplus yA(1)$. This step will be fulfilled separately for the two different algebras under consideration. Let us denote by $f_x$ and $f_y$ the degree $0$ right $A$-module generators of $xA(1)$ and $yA(1)$ correspondingly.

Suppose that $A=\kk\langle x,y\rangle/(xy,yx)$. Since $(xf_x)y=(yf_x)y=0$, we have $xf_x=q_xf_xx$ and $yf_x=q_yf_xx$ for some $q_x,q_y\in \kk$. We have shown before that $xM\cap yM=0$, i.e. $q_x=0$ or $q_y=0$. Analogously, $xf_y=p_xf_yy$ and $yf_y=p_yf_yy$, where one of the elements $p_x,p_y\in\kk$ is zero and another is nonzero. Moreover, it is clear that either $q_y=p_x=0$, $q_x\not=0$, $p_y\not=0$ or $q_x=p_y=0$, $q_y\not=0$, $p_x\not=0$. Then $M= \kk^m\oplus {}_{\gamma}(AxA\oplus AyA)$, where $\gamma$ is an automorphism of $A$ that sends $x$ to $q_xx+q_yy$ and $y$ to $p_xx+p_yy$.

If $q_y=p_x=0$, then ${}_{\gamma}(AxA)\ot_A{}_{\gamma}(AyA)\cong {}_{\gamma}(AyA)\ot_A{}_{\gamma}(AxA)\cong\kk$ and the definition of an $s$-homogeneous triple gives us isomorphisms $\alpha:\big({}_{\gamma}(AxA)\big)^{\ot_A^s}\rightarrow (AxA)$ and $\beta:\big({}_{\gamma}(AyA)\big)^{\ot_A^s}\rightarrow (AyA)$ that satisfy the equalities $\alpha(x^{\ot_A^s})*x=x*\alpha(x^{\ot_A^s})$ and $\beta(y^{\ot_A^s})*y=y*\beta(y^{\ot_A^s})$, where $*$ is the multiplication arising from the $A$-bimodule structure on ${}_{\gamma}(AxA\oplus AyA)$. Now it is clear that $\gamma=\Id_A$, $M=\kk^m\oplus AxA\oplus AyA$, and $\Lambda\cong\kk\langle x_1,\dots,x_m,y_1,y_2\rangle/(y_1^s,y_2^s)$. This algebra is $s$-Koszul, for example, by the results of \cite{MSV}, or \cite{GMMZ}.

If $q_x=p_y=0$, then analogous argument shows that $2\mid s$ and $q_y=\frac{1}{p_x}=q$ for some $q\in\kk^*$. It is not difficult to see that the graded algebra isomorphism from $A$ to $A$ that sends $x$ to $qx$ and $y$ to $y$ induces isomorphism of $s$-homogeneous triples $$(A,M,\vp)\cong \big(A,\kk^m\oplus {}_{\gamma_0}(AxA\oplus AyA),\vp_0\big),$$ where $\gamma_0$ is the automorphism of $A$ interchanging $x$ and $y$ and $\vp_0$ is the corresponding homomorphism from the definition of an $s$-homogeneous triple. Thus, $\Lambda\cong\kk\langle x_1,\dots,x_m,y_1,y_2\rangle/\big((y_1y_2)^t,(y_2y_1)^t\big)$ for $t=\frac{s}{2}$ is again $s$-Koszul.

Suppose now that $A=\kk\langle x,y\rangle/(x^2,y^2)$. Since $(xf_x)x=(yf_x)x=0$, we have $xf_x=q_xf_yx$ and $yf_x=q_yf_yx$ for some $q_x,q_y\in \kk$. Analogously to the previous case we have either $yf_x=xf_y=0$, $xf_x=q_xf_yx$, and $yf_y=p_yf_xy$ for some $q_x,p_y\in\kk^*$ or $xf_x=yf_y=0$, $yf_x=q_yf_yx$, and $xf_y=p_xf_xy$ for some $q_y,p_x\in\kk^*$. Then it is not difficult to show that $M= \kk^m\oplus {}_{\gamma}(A_{>0})$, where $\gamma$ is an automorphism of $A$ that sends $x$ to $q_xx+q_yy$ and $y$ to $p_xx+p_yy$.

If $q_y=p_x=0$, then $\big({}_{\gamma}(A_{>0})\big)^{\ot_A^s}\cong \kk^{2^s-2}\oplus L$, where $L$ is the $A$-bimodule generated by the elements $f_x^{\ot_A^s}$ and $f_y^{\ot_A^s}$.
It is easy to see that $x(l\Ann_{y}L)=0$, and hence $L\not\cong A_{>0}$.

If $q_x=p_y=0$, then $\big({}_{\gamma}(A_{>0})\big)^{\ot_A^s}\cong \kk^{2^s-2}\oplus L$, where $L$ is the $A$-bimodule generated by the elements
$$f_x^{\ot_A^{s-2\upper{\frac{s-1}{2}}}}\ot_A(f_y\ot_Af_x)^{\ot_A^{\upper{\frac{s-1}{2}}}}\mbox{ and }f_y^{\ot_A^{s-2\upper{\frac{s-1}{2}}}}\ot_A(f_x\ot_Af_y)^{\ot_A^{\upper{\frac{s-1}{2}}}}.$$
If $2\mid s$, then $L\not\cong A_{>0}$, because $x(l\Ann_{y}L)=0$.
If $2\nmid s$, then as before we get $q_y=\frac{1}{p_x}=q$ for some $q\in\kk^*$ and $(A,M,\vp)\cong \big(A,\kk^m\oplus {}_{\gamma_0}(A_{>0}),\vp_0\big)$, where $\gamma_0$ is the automorphism of $A$ interchanging $x$ and $y$ and $\vp_0$ is the corresponding homomorphism from the definition of an $s$-homogeneous triple. Thus, $\Lambda\cong\kk\langle x_1,\dots,x_m,y_1,y_2\rangle/\big((y_1y_2)^ty_1,y_2(y_1y_2)^t\big)$ for $t=\frac{s-1}{2}$ is $s$-Koszul.

It is clear that if $s\ge 3$, then $\Lambda^!$ is not $s$-Koszul in all cases. Thus, we obtain the following result.

\begin{theorem} Suppose that $\Lambda\in\HA(\kk,s)$.\\
1. If $(\Lambda^!)^{(s)}\cong\kk\langle x,y\rangle/(xy,yx)$, then either
$$\Lambda\cong\kk\langle x_1,\dots,x_m,y_1,y_2\rangle/(y_1^s,y_2^s)$$
for some $m\ge 0$ or $s=2t$ and
$$\Lambda\cong\kk\langle x_1,\dots,x_m,y_1,y_2\rangle/\big((y_1y_2)^t,(y_2y_1)^t\big)$$
for some $t\ge 1$ and $m\ge 0$.\\
2. If $(\Lambda^!)^{(s)}\cong\kk\langle x,y\rangle/(x^2,y^2)$, then $s=2t+1$ and
$$\Lambda\cong\kk\langle x_1,\dots,x_m,y_1,y_2\rangle/\big((y_1y_2)^ty_1,y_2(y_1y_2)^t\big)$$
for some $t\ge 1$ and $m\ge 0$.

In particular, if $(\Lambda^!)^{(s)}\in\left\{\kk\langle x,y\rangle/(xy,yx),\kk\langle x,y\rangle/(x^2,y^2)\right\}$, then $\Lambda$ is $s$-Koszul and $\Lambda^!$ is $s$-Koszul only in the case $s=2$.
\end{theorem}

\end{document}